\input amstex.tex
\documentstyle{amsppt}
\UseAMSsymbols
\TagsAsMath
\widestnumber\key{ASSSS}
\magnification=\magstephalf
\pagewidth{6.2in}
\vsize8.0in
\parindent=6mm
\parskip=3pt
\baselineskip=16pt
\tolerance=10000
\hbadness=500
\NoRunningHeads
\loadbold

\document
\def \rectangle#1#2{\hbox{\vrule\vbox to #2
              {\hrule\hbox to #1{\hfil}\vfil\hrule}\vrule}}
\def\sq{\,\,\rectangle{7pt}{7 pt}\,\,}

\def\ring{\!\circ\!}
\NoBlackBoxes
\nologo
\topmatter
\title A remark on Global existence for small initial data\\
of the minimal surface equation in Minkowskian space time
\endtitle
\author Hans Lindblad \endauthor
\address 
University of California at San Diego
\endaddress 
\endtopmatter

\head 1.Introduction\endhead 
We show that the nonlinear wave equation corresponding 
to the minimal surface equation in Minkowski space time 
$$
\frac{\partial}{\partial t} 
\frac{\phi_t}{\sqrt{1+|\nabla_x \phi|^2-\phi_t^2}} 
-\sum_{i=1}^n\frac{\partial}{\partial x^i} 
\frac{\phi_i}{\sqrt{1+|\nabla_x \phi|^2-\phi_t^2}}=0\tag 1.1  
$$
where $\phi_i=\partial \phi/\partial x^i$, $\phi_t=\partial \phi/\partial t$,
has global solutions for sufficiently small initial data:
$$
\phi\big|_{t=0}=\varepsilon f,\qquad 
 \phi_t\big|_{t=0}=\varepsilon g \tag 1.2 
$$
$f\in C_0^\infty$ and $g\in C_0^\infty$, i.e. (1.1)-(1.2) has for fixed $f$ and
$g$ a solution for all $t\geq 0$ if $\varepsilon>0$ is sufficiently small. 
This is an interesting model in Lorentzian 
geometry proposed to me by Hamilton\cite{Ha1}. 
It is also the equation for a membrane, 
in field theory, see Hoppe\cite{Ho1}. 
Also Huisken and Struwe\cite{HS1} have some 
recent results related to local existence for (1.1). 

What makes the proof go through also in the physically 
important case of two space dimensions ($\phi$ itself 
corresponds to the third space dimension) is that
the nonlinear terms satisfies the so called 
"null condition" of Christodoulou\cite{C1} and 
Klainerman\cite{K2,K4}. 
The purpose of this note is to present two simple proofs making use of the 
extra symmetries of the equation. 
The first proof uses a version of the method of \cite{K2}, that works also in two 
space dimensions. For equations in divergence form
we can get a good $L^2$ estimate for the solution itself, 
see \cite{L1}, that replaces the 
conformal energy estimate used in \cite{K2}. 
The second proof uses a simplified version \cite{C2} of the method of \cite{C1} and
it works also in one space dimension due to 
that the equation satisfies a "double null condition", see (1.5).
The first proof does not work in the case of one space dimension but it has
the advantage that it does not require compact support of initial data
but merely some decay at infinity. 

Let 
$$
\square=\partial_t^2-\sum_{i=1}^3\partial_i^2,\qquad
\text{and}\qquad
Q_{00}(\phi,\psi)=\phi_t \psi_t-\sum_{i=1}^n \phi_i\psi_i\tag 1.3 
$$
be a null form. We can write (1.1) as a 
wave equation with a right hand side in divergence form 
$$
\square \phi=\partial_t \big(\phi_t F(Q_{00}(\phi,\phi))\big)
-\sum_{i=1}^n \partial_i \big(\phi_i F(Q_{00}(\phi,\phi))\big),\tag 1.4
$$
where $F(Q)=-1+{1}/{\sqrt{1-Q}}$.
We can alternatively write it with null forms 
$$
\square \phi= -\frac{Q_{00}\big(\phi,
Q_{00}(\phi,\phi)\big)}{2(1-Q_{00}(\phi,\phi))}\tag 1.5 
$$
The idea behind the proofs is that we expect solutions of the nonlinear 
wave equation to decay like solutions of the linear homogeneous wave equation
$\square \phi=0$, i.e. 
$|\partial^\alpha \phi|\leq C\varepsilon (1+t)^{-(n-1)/2}$. 
This will make the right hand sides of the nonlinear wave equations
small for large $t$ and hence these equations will be close to the homogeneous case
so we can close the argument. 
What makes the proof go through in the lower dimensional case is that 
for a null form there is an additional cancellation which 
leads to an additional decay of a factor of $(1+t)^{-1}$:
$|\partial^\alpha Q(\phi,\phi)|\leq C\varepsilon^2(1+t)^{-(n-1)-1}$.

\head 2. The proof using vector fields\endhead 
The idea of this argument is to exploit that solutions of linear wave equations
$\square v=0$ satisfies the decay estimate, see e.g. \cite{H\"o2}, 
$$
|v(t,x)|\leq C(f,g)\varepsilon (1+t+|x|)^{-(n-1)/2}(1+|t-|x||)^{-(n-1)/2}, \tag 2.1 
$$
where $C(f,g)$ is a constant depending on some weighted Sobolev norm of initial data
 $(f,g)$. We will use vector fields to obtain this kind
of decay. 
For $(t,x)\in \bold{R}^{1+n}$ denote $\partial_t$ by $\partial_0$
and $\partial_{x_j}$ by $\partial_j$ for $j=1,...,n$. Let 
 $$
\Gamma_{jk}={\lambda}_jx_j\partial_k-{\lambda}_kx_k\partial_j,\quad\text{for}\quad 
i\neq j,\qquad\text{and}
\qquad \Gamma_{00}=\sum_0^n{x_j\partial_j}\tag 2.2
 $$
where ${\lambda}=(1,-1,...,-1)$ and $x_0=t$.
$\Gamma_{jk}$, for $(j,k)\neq (0,0)$
 are the vector fields associated with the 
Lorentz group that all commute with $\square$ and $\Gamma_{00}$
is the scaling vector fields whose commutator is 
$[\Gamma_{00} ,\square]=-2\square$. 
$\Gamma$ will symbolically stand for any of the vector fields 
$\Gamma_{ij}$ or $\partial_j\,$, $i,j=0,1,...,n$ and we will write 
$\Gamma^I$ for a product of $|I|$ of such vector fields. 
We note that $[\partial_i,\Gamma_{jk}]$ is either $0$ or else 
equal to $\pm\partial_l$ for some $l$.
The operators $\{\Gamma_{jk}\}$ span the tangent space at every 
point where $t\neq |x|$. But when $t=|x|$ they only span the tangent
space of the cone $t=|x|$, so we have, see \cite{K1,2}, \cite{H\"o1,2} or
\cite{L1}; 
 $$
|\partial \phi|\leq C(|t-|x||)^{-1}\sum_{}{|\Gamma_{ij} \phi|},
\qquad 
|\partial \phi|\leq C(t+|x|)^{-1}\big(|x||\partial_t
\phi|+\sum_{}{|\Gamma_{ij}\phi|}\big),\tag 2.3
$$
where  $|\partial w |^2=\sum_0^n{|\partial_j w|^2}$.
Now we also need to calculate the commutator of $\Gamma$ with a 
null form $Q$. However in these commutators other null forms than 
(1.3) will come up. For $0\leq i,j\leq n$ and $i\neq j$ let 
$$
Q_{ij}(\phi,\psi)=(\partial_i \phi)\partial_j\psi-
(\partial_j \phi)\partial_i\psi,\tag 2.4
$$
Let $Q$ symbolically stand for any of the null forms 
(1.3) or (2.4). Then
$$
\Gamma Q(\phi,\psi)= Q(\Gamma \phi,\psi)+Q(\phi,\Gamma\psi)
+\sum a_{ij} Q_{ij}(\phi,\psi) \tag 2.5
$$
for some constants $a_{ij}$. 
For a null form we have better decay close to the light cone, see\cite{K1,2}
$$
|Q(\phi,\psi)|\leq 
C(1+|t|+|x|)^{-1} \big(|\partial \phi| |\Gamma \psi|+|\Gamma \phi| |\partial \psi |
\big) \tag 2.6
$$
where here $|\Gamma \phi|^2=\sum_{|I|=1} |\Gamma^I \phi|^2$. 
The decay of $\phi$ will be obtained from that $\Gamma^I\phi$ also satisfies 
nonlinear wave equations of the same form. 
Applying $\Gamma^I$ to (1.4) we get 
$$
\square \Gamma^I\phi=\sum_{j=0}^n  \partial_j 
\bigg( \sum_{\Sb k\geq 3,\, |I_1|+...+|I_k|\leq
|I|\\ |I_i|\leq |I|/2 , \,i<k,\, \,|I_k|\leq |I|\endSb} F_{jk i_1...i_k I_1...I_k} (Q_{00}(\phi, \phi))
(\partial_{i_1}\Gamma^{I_1}\phi)(\partial_{i_2}\Gamma^{I_2}\phi)\cdot\cdot\cdot
(\partial_{i_k}\Gamma^{I_k}\phi)\bigg) \tag 2.7
$$
Similarly, from (1.5) we get 
$$\multline 
\square \Gamma^I \phi+\frac{Q_{00}\big(\phi,
Q_{00}(\phi,\Gamma^I\phi)\big)}{1+Q_{00}(\phi,\phi)}\\
=\sum_{\Sb k\geq 3,\, |I_1|+...+|I_k|\leq
|I|+1,\\ |I_i|\leq (|I|+1)/2,\, i<k,\,\, |I_k|\leq |I|\endSb } 
G_{k i_1...i_k I_1...I_k} (Q_{00}(\phi, \phi))
(\partial_{i_1}\Gamma^{I_1}\phi)(\partial_{i_2}\Gamma^{I_2}\phi)\cdot\cdot\cdot
(\partial_{i_k}\Gamma^{I_k}\phi)
\endmultline 
\tag 2.8
$$
and
$$
\multline 
\square \Gamma^I \phi \\
=\!\!\!\!\!\!\!\!\!\!\!\!\!\!\!\!\!\!\!\!\sum_{\Sb\,\,\,\, \,\,\,\,\, \,\,\,\,\,k\geq 3,\, |I_1|+...+|I_k|\leq |I|\endSb }
\!\!\!\!\!\!\!\!\!\!\!\!\!\!\!\!\!\!\!
 H_{I_1...I_k}^{k i_0...i_k}(Q_{00}(\phi,\phi))    Q_{i_0 i_1}\Big(\Gamma^{I_1}\phi,
Q_{i_2j_3}(\Gamma^{I_2}\phi,\Gamma^{I_3}\phi)\Big)
Q_{i_4 i_5}(\Gamma^{I_4} \phi,\Gamma^{I_5} \phi) \cdot\cdot \cdot 
Q_{i_{k-1} i_k} (\Gamma^{I_{k-1}} \phi,\Gamma^{I_k} \phi) 
\endmultline \tag 2.9
$$
where $k$ is odd and if $k=3$ then this is to be interpreted as
that the factor $Q_{i_4 i_5}(\Gamma^{I_4} \phi,\Gamma^{I_5} \phi)$
is absent. 
	
The proof will use the energy inequality applied to derivatives of
the solution and some decay estimates. 
The energy inequality, see \cite{K3} or \cite{H\"o2},
for a solution of 
 $$
\sq w+\sum_{j,k=0}^n{{\gamma}^{jk}(t,x)\partial_j\partial_k w}=F, 
\qquad |{\gamma}|=\sum{|{\gamma}^{jk}|}\leq \frac{1}{2}\tag 2.10
 $$
says that 
 $$
\|\partial w(t,\cdot)\|_{L^2}\leq 2
\exp{\tsize{(\int_0^t{2|{\gamma}^{\prime}(\tau)|\,d\tau}})}
\|\partial w(0,\cdot)||_{L^2}
+2\int_0^t \exp{\tsize{(\int_s^t{2|{\gamma}^{\prime}(\tau)|\,d\tau}})}
{||F(s,\cdot)||_{L^2}\,ds}, 
\tag 2.11
 $$
where $|{\gamma}^{\prime}(t)|=\sum_{i,j,k}\sup{|\partial_i{\gamma}^{jk}(t,\cdot)|}$.
It was not stated exactly like this in \cite{H\"o2} but (2.11)
follows from the proof of the version there. 
As a consequence of the energy inequality we also have the following 
estimate: If 
 $$
\sq w=\sum_0^n \partial_j F_j, \qquad 
w\big|_{t=0}=\varepsilon f,\quad w_t\big|_{t=0}=\varepsilon g\tag 2.12
 $$
then 
 $$
||w(t,\cdot)||_{L^2}\leq
\sum_{j=0}^n\int_0^t{||F_j(s,\cdot)||_{L^2}\,ds}+C(f,g,F_0(0,\cdot))m(t)\varepsilon\tag 2.13
$$
where $m(t)=1$ if $n\geq 3$, $m(t)=\log{(2+t)}$ if $n=2$ and
 $C(f,g,F_0(0,\cdot))$ stands for some constant depending on some weighted Sobolev
norm of initial data $f$ and $g$.
The proof is a trick used in \cite{L1}; if
$\square v_j=F_j$ then $\square(w-\sum\partial_j v_j)=0$
so $\|w(t,\cdot)\|_{L^2}$ is bounded by 
 $\sum\|\partial v_j(t,\cdot)\|_{L^2}$, which can be estimate
by the energy inequality (2.11), 
plus the norm for a solution of a linear homogeneous equation, 
which can be obtained from e.g. (2.1). 
We will also need an $L^1-L^\infty$ estimate of
H\"ormander\cite{H\"o1}, (see also Klainerman\cite{K1,2} for 
an earlier version and \cite{L1} for a simple proof): The solution $w$ of 
 $$
\sq w=F,\qquad   w\big|_{t=0}=\varepsilon f,
\quad w_t\big|_{t=0}=\varepsilon g\tag 2.14
 $$
satisfies 
 $$
|w(t,x)| \leq C(1+t+|x|)^{-(n-1)/2}\Big(\sum_{|I|\leq n-1}{\int_0^t
{||(\Gamma^I F)(s,\cdot)/(1+s+|\cdot|)^{(n-1)/2}||_{L^1}\,ds}}
+ C(f,g)\varepsilon \Big)
\tag 2.15
 $$
Here the estimate for the 
linear homogeneous part, the second term, is (2.1). 
Whereas the proof of the energy inequality is merely integration by parts the 
proofs of the decay estimates (2.15) and (2.1) requires a detailed analysis
of the fundamental solution or stationary phase.

Let $N\geq 2n+1$ and $\delta=0$, if $n\geq 3$, and $0<\delta<1/2$ fixed,if $n=2$.
We will now prove that  
$$
\align 
M_1(t)&=\sum_{|I|\leq N}\|\partial\Gamma^I\phi(t,\cdot)\|_{L_2}\leq 
K\varepsilon (1+t)^{\delta},\tag 2.16\\
M_2(t)&=\sum_{|I|\leq N}\|\Gamma^I\phi(t,\cdot)\|_{L_2}
\leq K \varepsilon (1+t)^{\delta},\tag 2.17\\
N_1(t)&=\sum_{|J|\leq (N+1)/2}
\|\partial\Gamma^J\phi(t,\cdot)\|_{L^\infty}
\leq K\varepsilon(1+t)^{-(n-1)/2},\tag 2.18\\
N_2(t)&=\sum_{|J|\leq (N+1)/2+1} \|\Gamma^J\phi(t,\cdot)\|_{L^\infty}
\leq K\varepsilon(1+t)^{-(n-1)/2},\tag 2.19
\endalign 
$$
if $K$ is sufficiently large and $\varepsilon $ is sufficiently small.
We observe that the bound for $N_1$ is a consequence of the
bound for $N_2$ since in particular we can take one 
factor of $\Gamma=\partial$. 
(2.11) applied to (2.8) gives
$$
M_1(t)\leq 
C\varepsilon \exp{\tsize{(\int_0^t{N_1(\tau)^2\,d\tau}})}
+ \int_0^t \exp{\tsize{(\int_s^t{N_1(\tau)^2\,d\tau}})}
C(N_1(s)) N_1(s)^2 M_1(s)\, ds
\tag 2.20
$$
and (2.13) applied to (2.7) gives 
$$
M_2(t)\leq 
C\Big( C(f,g)m(t)\varepsilon + \int_0^t C(N_1(s)) N_1(s)^2 M_1(s)\, ds\Big)
\tag 2.21
$$
Finally, (2.15) applied to (2.9) using (2.6) and Cauchy
Schwartz inequality gives 
$$
N_2(t)\leq 
C(1+t)^{-(n-1)/2}\Big( C(f,g)\varepsilon +
\int_0^t \frac{(N_1(s)+N_2(s))}{(1+s)^{(n-1)/2+1}} (M_1(s)+M_2(s))^2 
\,
ds\Big)
\tag 2.22
$$
if $(N+1)/2+1+n-1\leq N$, i.e. $N\geq 2n+1$.  
What will make the argument work also for $n=2$
is that the null condition gave an extra power of $(1+s)^{-1}$ in the integral (2.22).
In fact because we have a double null condition we 
actually have one more power but we have no use of this here.
The rest of the argument is now by continuity.
We know that (2.16)-(2.19) are true for $t=0$ 
if $K$ is large enough and we know from 
the local existence theorem for hyperbolic equations, 
see e.g. \cite{H\"o2}, that 
these quantities are continuous as long as they are bounded.
We now assume that $T_1$ is the largest number such that (2.16)-(2.19) 
are true for $t\leq T_1$ 
and show that these bounds together with (2.20)-(2.22)
implies stronger bounds if $K$ is sufficiently large and $\varepsilon$ is
sufficiently small. Hence by continuity we conclude that the bounds 
(2.16)-(2.19) must hold for $t\leq T_2$ where $T_2>T_1$, contradicting 
the maximality of $T_1$. 
To simplify notation, let us now only deal
with the most sensitive case $n=2$.
Then
$$
\exp{\big(\tsize{\int_s^t N_1(\tau)^2\, d\tau}\big)}
\leq \exp{\big(\tsize{K^2\varepsilon^2 \int_s^t(1+\tau)^{-1}\, d\tau }\big)}=
\exp{\big(K^2\varepsilon^2\ln{\big(\frac{1+t}{1+s}\big)}\big)}
=\left(\frac{1+t}{1+s}\right)^{K^2\varepsilon^2}\tag 2.23
$$
so it follows from (2.20): 
$$
M_1(t)\leq C\varepsilon (1+t)^{K^2\varepsilon^2}
+\int_0^t C\varepsilon^3 K^3
\left(\frac{1+t}{1+s}\right)^{K^2\varepsilon^2}
\!\!\!\!\!(1+s)^{\delta-1}\, ds\leq 
K\varepsilon(1+t)^\delta/2 \tag 2.24
$$
if $K$ is sufficiently large and $\varepsilon$ is sufficiently small. 
Similarly, from (2.21) we get 
$$
M_2(t)\leq C\varepsilon\log{(2+t)} +\int_0^t K^3\varepsilon^3 (1+s)^{\delta-1}\, d s
\leq K\varepsilon(1+t)^\delta /2\tag 2.25
$$
if $K$ is sufficiently large and $\varepsilon$ is sufficiently small. Finally from
(2.22) we get 
$$
N_2(t)\leq C(1+t)^{-1/2}\Big(\varepsilon + 
\int_0^t K^3\varepsilon^3(1+s)^{2\delta-2}\, ds\Big) \leq 
K\varepsilon(1+t)^{-1/2}/2\tag 2.26
$$
if $K$ is sufficiently large and
$\varepsilon$ is sufficiently small since $0<\delta<1/2$. 
This concludes the proof.

\def\ka{\kappa\,}
\head 3. The proof using conformal inversion\endhead
We will reduce the global problem to a local problem, for which small data existence
is known, using a conformal inversion or Kelvin transform.
Let 
$\kappa:\bold{R}^{1+n}\ni (s,y)\to (t,x)\in \bold{R}^{1+n}$ and
$\tilde{\phi}$ be defined by 
$$
\tilde{\phi}=\phi\ring\ka \rho^{-\alpha},
\qquad x^i=\kappa^i(s,y)=y^i/\rho,\qquad \rho=s^2-|y|^2,
\qquad \alpha=\tfrac{n-1}{2},\tag 3.1
$$
where $x^0=t$ and $y^0=s$. 
Let $m_{ij}=m^{ij}$ be the Minkowski metric $m_{00}=1$,
$m_{ii}=-1$ if $j\geq 1$ and $m_{ij}=0$ if  $i\neq j$.
Then 
$\partial_i \kappa^k=\delta_i^{\,\,\, k}/\rho-2y^j y^k m_{ij}/\rho^2$.
If the metric in the $x$ coordinates is $m_{ij}$ then the pull-back
metric in the $y$ coordinates is given by 
$$
g_{ij}=m_{kl}(\partial_i \kappa^k )\partial_j \kappa^l=m_{ij}/\rho^2,
\qquad g^{ij}=m^{ij}\rho^2\tag 3.2
$$
Then the norm is invariant:
$
g^{ij}(\partial_i \phi\ring\ka )\partial_j \psi\ring \ka=
m^{ij}\phi_i\ring\ka \psi_j\ring\ka
$
where $\partial_i \phi\ring\ka=\partial\phi\ring\ka/\partial y^i$
and $\phi_i=\partial \phi/\partial x^i$,  
i.e. if $Q_{00}$ is the null form (1.3) then 
$$
Q_{00}(\phi,\psi)\ring\ka=\rho^2 Q_{00}(\phi\ring\ka,\psi\ring\ka).\tag 3.3
$$
Expressing $\square$ in the $y$ coordinates we get 
$$
(\square \phi)\ring \kappa=\square_g (\phi\ring\kappa)=
(\det{g})^{-1/2}\partial_i \big( g^{ij}(\det{g})^{1/2} \partial_j
\phi\ring\kappa\big)=\rho^{n+1}m^{ij}\partial_i
 \big(\rho^{-(n-1)} \partial_j\phi\ring\kappa\big)
\tag 3.4
$$
Since the operator $\square_g-(n-1)R/4n$, where $R$ is the scalar curvature,
is conformally covariant and the scalar curvature in both cases vanishes, we obtain 
$$
\square\tilde{\phi}= \rho^{-\alpha-2}(\square \phi)\ring\kappa \tag 3.5 
$$
This is also follows from (3.4) by 
an easy calculation using that $\square \rho^{-\alpha}=0$.
Hence by (1.5), (3.5) and (3.3) 
$$
\square \tilde{\phi}= -\rho^{-\alpha} 
\frac{Q_{00}\big(\rho^\alpha\tilde{\phi},
 \rho^2 Q_{00}(\rho^\alpha\tilde{\phi},\rho^\alpha\tilde{\phi})\big)}{
1-\rho^2 Q_{00}(\rho^\alpha\tilde{\phi},\rho^\alpha\tilde{\phi})}\tag 3.6
$$
Since $m^{ij} (\partial_i \tilde{\psi}) \partial_j\rho=2(\Gamma_{00} \tilde{\psi})$,
where $\Gamma_{00}=s\partial_s +y^i\partial_i$, 
and $\Gamma_{00}\rho=2\rho$ we obtain 
$$
Q_{00}(\rho^\beta\tilde{\phi},\rho^\gamma \tilde{\psi})=
\rho^{\beta+\gamma-1}\Big(
\rho Q_{00}(\tilde{\phi},\tilde{\psi})+2\beta
\tilde{\phi}(\Gamma_{00}\tilde{\psi})
+2\gamma\tilde{\psi} (\Gamma_{00} \tilde{\phi})
+4\gamma\beta\tilde{\phi}\tilde{\psi} \Big),\tag 3.7
$$
Furthermore, $\Gamma_{00}Q_{00}(\tilde{\phi},\tilde{\phi})=
2Q_{00}(\tilde{\phi},\Gamma_{00}\tilde{\phi})-
2Q_{00}(\tilde{\phi},\tilde{\phi})$. 
Using (3.7) twice we see that we can factor out 
$\rho^{3\alpha}$ from ${Q_{00}\big(\rho^\alpha\tilde{\phi},
 \rho^2 Q_{00}(\rho^\alpha\tilde{\phi},\rho^\alpha\tilde{\phi})\big)}$ so (3.6)
has the general form: 
$$\multline 
\square \tilde{\phi}
=-\rho^{2\alpha}\frac{\rho^2
Q_{00}\big(\tilde{\phi},Q_{00}(\tilde{\phi},\tilde{\phi})\big) 
+\rho c_1\tilde{\phi}Q_{00}(\tilde{\phi},\Gamma_{00}\tilde{\phi})
+\rho c_2 \tilde{\phi}^2 \Gamma_{00}^2\tilde{\phi}}
{1-\rho^{2\alpha}\big(\rho^2 Q_{00}(\tilde{\phi},\tilde{\phi})
+4\alpha\rho (\alpha\tilde{\phi}^2+\tilde{\phi}\Gamma_{00}\tilde{\phi})\big)}\\
+\rho^{2\alpha}\frac{\rho(c_3\tilde{\phi}+c_4 \Gamma_{00}\tilde{\phi})
Q_{00}(\tilde{\phi},\tilde{\phi})+\tilde{\phi}(c_5 \tilde{\phi}^2+
c_6\tilde{\phi}\Gamma_{00}\tilde{\phi}+c_7(\Gamma_{00}\tilde{\phi})^2)}
{1-\rho^{2\alpha}\big(\rho^2 Q_{00}(\tilde{\phi},\tilde{\phi})
+4\alpha\rho (\alpha\tilde{\phi}^2+\tilde{\phi}\Gamma_{00}
\tilde{\phi})\big)}
\endmultline \tag 3.8
$$
for some constants $c_1,...,c_7$. We can write this as 
$$
\square \tilde{\phi}=\rho^{n-1}
\big( \rho B^{ij}(s,y,\tilde{\phi},\partial\tilde{\phi})\partial_i \partial_j \tilde{\phi}
+T(s,y,\tilde{\phi},\partial\tilde{\phi})\big),\qquad \rho=s^2-|y|^2, \tag 3.9
$$
where $B^{ij}$ and $T$ are smooth functions of $(s,y)$,
$\tilde{\phi}$ and $\partial_i \tilde{\phi}$, for $i=0,...,n$
vanishing to second respectively third order at
$(\tilde{\phi},\partial\tilde{\phi})=(0,0)$;
$|B^{ij}(s,y,\tilde{\phi},\partial\tilde{\phi})|
\leq C (|\tilde{\phi}|+|\partial\tilde{\phi}|)^2$ 
and $|T(s,y,\tilde{\phi},\partial\tilde{\phi})|
\leq C (|\tilde{\phi}|+|\partial\tilde{\phi}|)^3$.
The importance of (3.9) is that it is nonsingular at $\rho=0$. 

Having derived the transformation of the equation let us
explain in what region it is applied. The transformation 
$\kappa$ maps the interior of the forward light cone 
$\{(s,y);\,s>|y|\}$ onto itself $\{(t,x);\,t>|x|\}$ and the boundary
$\{(s,y);\,s=|y|\}$ to infinity. 
Its inverse $\kappa^{-1}(t,x)=\big(t/(t^2-|x|^2),x/(t^2-|x|^2)\big)=(s,y)$
maps forward light cones $\Lambda_c=\{(t,x);\, t-|x|\geq c>0\}$
into backward light cones 
$\tilde{\Lambda}_c=\{(s,y);\, s+|y|\leq 1/c,\, s-|y|\geq 0\}$
intersected with $\{(s,y);\, s\geq |y|\}$. 
The inverse is smooth on $\Lambda_c$ since $t^2-|x|^2\geq c^2>0$ there.
Note also that hyper planes $\tilde{\Cal H}_b=\{(s,y);\, s=1/2b,\, |y|\leq 1/2b\}$ are 
transformed to hyperboloids  ${\Cal H}_b=\{(t,x);\, (t-b)^2-|x|^2=b^2\}$.
We will now use this transformation in  $\Lambda_c$ 
which by the inverse is mapped to the compact region $\tilde{\Lambda}_c$
where we can use the standard local existence theorem for (3.9) to, after
transforming back, obtain a 
global solution of our original equation (1.5) in $\Lambda_c$. 
Let us now explain how we 
can reduce it to a problem in $\Lambda_c$.
First by scaling, $\phi_a(t,x)=\phi(at,ax)/a$ 
is a solution of (1.5) if $\phi$ is, so we may assume that data 
(1.2) are supported in the set $\{x;|x|\leq 1\}$. 
Secondly, we can translate the solution in the 
time direction so initial conditions are attained when $t=a>1$.
By the local existence theorem we have, 
if initial data when $t=a$ are sufficiently small,  a solution $\phi$ to (1.5), 
for $|t-a|\leq a$, and it 
is as small as we wish there. Furthermore, by Huygens' principle the solution
vanishes outside the forward light cone $\Lambda_{a-1}$, for $t\geq a$. 
We now want to show that the solution $\phi$ of (1.5) extends
to a solution in all of $\Lambda_{a-1}$ for $t\geq a$ 
by showing that we have a local solution 
$\tilde{\phi}$ of (3.9) in the compact set $\tilde{\Lambda}_{a-1}$.
We will now describe how to obtain initial conditions for (3.9). 
We pick a particular hyperboloid 
 ${\Cal H}_b$, where $b=(a-1/a)/2>0$, 
that intersects with the plane $t=a$ exactly when $|x|=1$, and that is 
transformed by the inverse of $\kappa$ to the plane 
$\tilde{\Cal H}_b$. The intersection of ${\Cal H}_b$ with the support 
of $\phi$, ${\Cal H}^1_b=\{(t,x)\in {\Cal H}_b;\, |x|\!\leq\! 1\}$, 
is contained it the forward light cone $\Lambda_{a-1}$ intersected with the set 
where $t\leq a$. Hence we have a smooth solution  
$\phi$ on ${\Cal H}_b^1$ which is as small as we wish there and which vanishes 
on ${\Cal H}_b\setminus {\Cal H}_b^1$. It follows that $\phi$ and its derivatives 
restricted to ${\Cal H}_b$ are transformed onto smooth initial conditions for
$\tilde{\phi}$ on $\tilde{\Cal H}_b$, which are as small as we wish. 
The local existence theorem with these initial conditions on $\tilde{\Cal H}_b$
gives us a smooth solution $\tilde{\phi}$ of (3.9) in
$\{(s,y); s\!\leq\! 1/2b , \, s\!\geq\! |y|\}$ and hence in 
 $\tilde{\Lambda}_{a-\!1}$ for $s\!\leq \! 1/2b$, if initial conditions on 
$\tilde{{\Cal H}}_b$ are sufficiently small. Transforming back gives us a global 
smooth solution $\phi$ of (1.5) in $\Lambda_{a-\!1}$ and hence for all $t\!\geq \!a$. 
This concludes the proof.

\subheading{Acknowledgments} I would like to thank Demetri Christodoulou,
Richard Hamilton and Jim Isenberg for stimulating discussions.


\Refs\ref\no [CC]\by 
Y. Choquet-Bruhat and D. Christodoulou\paper  
Existence of global solutions of the Yang Mills, Higgs and spinor field equations 
in $3+1$ dimensions\jour 
Ann. Sci. École Norm. Sup. \vol (4) 14 \yr 1981 \pages 481--506 
\endref 
\ref \no [C1]\by Christodoulou\paper global solutions of nonlinear 
hyperbolic equations for small initial data\jour
Comm. Pure Appl. Math.\vol 39\yr 1986\pages 267-282
\endref
\ref \no [C2]\bysame\paper Oral communication\yr 1999
\endref
\ref \no [C3]\bysame\paper Solutions globales des equations de champ de
Yang Mills\jour C.R. Acad. Sci. Paris (Series A)\vol 293 \yr 1981 \pages 139-141
\endref
\ref \no [CK1]\by D. Christodoulou and S. Klainerman
\book The nonlinear stability of Minkowski space-time\publ Princeton Univ. Press\endref
\ref\no [Ha1]\by R. Hamilton\paper Oral Communication\yr Oberwolfach 1994\endref 
\ref\no [Ho1]\by J. Hoppe\paper 
Some classical solutions of relativistic membrane equations in 
$4$-space-time dimensions.
\jour Phys. Lett. B 329 \yr 1994 \vol 1\pages  10-14\endref 
\ref\no [HS1]\by G. Huisken and M. Struwe\paper Oral
communication\yr 1999\endref
\ref \no [H\"o1]\by L. H\"ormander \book $L^1$, $L^\infty $ estimates for the
wave operator (Analyse math. et appl.)\publ
Gauthier-Villars\publaddr Paris 
\pages 211-234\yr 1988 \endref
\ref \no [H\"o2]\bysame \book Lectures on Nonlinear hyperbolic differential
equations\publ Springer Verlag \yr 1997\endref
\ref \no [JK1]\by F.\,John and S.\,Klainerman \paper Almost global existence to
nonlinear wave equations in three space dimensions
\jour Comm. Pure Appl. Math \vol 37\yr 1984\pages 443-455\endref
\ref \manyby S.\,Klainerman\no [K1]\paper Uniform decay estimates and the Lorentz
invariance of the wave equation\jour Comm. Pure Appl. Math\vol 38\yr
1985\pages 321-332\endref
\ref\bysame\no [K2]\paper The null condition and global
existence to nonlinear wave equations\jour Lectures in Applied
Mathematics \vol 23\yr 1986 \pages 293-326\endref
\ref \bysame\no [K3]\paper Global
existence for nonlinear wave equations\jour Comm. Pure Appl.
Math. \vol 33\yr 1980 \pages 43-101\endref
\ref\no [K4]\bysame \paper Long time behaviour of solutions to nonlinear wave equations
{\it Proceedings of the International Congress of Mathematics, 
(Warsaw, 1983)} \pages 1209-15
\publ PWN, Warsaw\yr 1984\endref 
\ref \no [LZ1]\by T. Li and Y. Zhou\paper Life-span of classical solutions 
to nonlinear wave equations in two space dimensions
\jour J. Math. Pure et Appl. \vol 73(3) \yr 1994 \pages
223--249 \endref 
\ref \no [LZ2]\by T. Li and Y. Zhou\paper Life-span of classical solutions 
to fully nonlinear wave equations in two space dimensions II
\jour J. Partial Diff. Eqs. \vol 6(1)\yr 1993\pages 17-38\endref 
\ref\no [L1]\by H. Lindblad\paper On the lifespan of solutions of nonlinear
wave equations with small initial data\jour Comm. Pure Appl. Math\vol 43\yr 1990\pages 445--472\endref
\ref\no [L2]\bysame\paper Global existence for nonlinear wave equations
\jour Comm. Pure Appl. Math \vol 45 (9)\yr 1992 \pages 1063-1096\endref
\endRefs
\enddocument